\documentclass[12pt]{amsart}
\usepackage{geometry}
\usepackage{graphicx}
\usepackage{amscd}
\usepackage{amssymb}
\usepackage{amsmath}
\usepackage[latin1]{inputenc}

\theoremstyle{plain}                                       %
\newtheorem{thm}{\quad Theorem}                            %
\newtheorem{cor}[thm]{\quad Corollary}                     %
\newtheorem{prop}[thm]{\quad Proposition}                  %
\theoremstyle{definition}                                  %
\newtheorem{alg}[thm]{\quad Algorithm}
\newtheorem{que}[thm]{\quad Question}                     %

\newcommand{\N}{{\Bbb N}}
\newcommand{\K}{{\Bbb K}}

\newcommand{\Z}{{\Bbb Z}}

\newcommand{\al}{{\alpha}}

\newcommand{\om}{{\omega}}

\begin{document}
{\title{Iterative processes related to Riordan arrays: The
reciprocation and the inversion of power series}
\author{Ana Luz\'{o}n }
\address{Departamento de Matem\'{a}tica Aplicada a los Recursos Naturales. E.T. Superior de
Ingenieros de Montes. Universidad Polit\'{e}cnica de Madrid.
28040-Madrid, SPAIN}\email{anamaria.luzon@upm.es}
\begin{abstract} We point out how Banach Fixed Point Theorem, and
the Picard successive approximation methods induced by it, allows
us to treat some mathematical methods in Combinatorics. In
particular we get, by this way, a proof and an iterative algorithm
for the Lagrange Inversion Formula.
\end{abstract}
\maketitle

\section{Introduction: beginning with a simple question.} \label{S:intro}

\noindent The results in this paper are consequences of special
interpretations, as fixed point problems, of the two classical
reversion processes in the realm of formal power series: the {\it
reciprocation}, i. e. the reversion for the Cauchy product, and
the {\it inversion}, i. e. the reversion for the composition of
series.

\noindent The case of the reciprocation was studied in \cite{teo}
and \cite{BanPas}. To unify our approach we also survey herein
some of the previous results on this topic.


\noindent The aim of this paper is to show how the (Picard)
successive approximation method induced by Banach's Fixed Point
Theorem (and some mild generalization) allows us to treat some
mathematical methods in combinatorics getting so some associated
algorithms. In particular:


1) To construct all elements in the Riordan group as a consequence
of the iterative process obtained to calculate the reciprocal  of
any power series admitting it. Presenting also an algorithm and
one pseudo-code description for it.

2) To construct, approximatively, the inverse of any power series
(admitting it) in such a way that the Lagrange Inversion Formula
can be first predicted and finally proved. We also describe the
corresponding algoritm.

\noindent For completeness we are going to recall the metric fixed
point theorems we will use, see, for example, \cite{DG82} for the
first one and   \cite{Sotomayor} in page 212 for the second one.

\noindent {\bf Banach Fixed Point Theorem (BFPT).} Let $(X,d)$ be
a complete metric space and $f:X\rightarrow X$ contractive. Then
$f$ has a unique fixed point  $x_0$ and $f^n(x)\rightarrow x_0$
for every $x\in X$.

\noindent {\bf Generalized Banach Fixed Point Theorem (GBFPT).}
Let $(X,d)$  be a complete metric space. Suppose
$\{f_n\}_{n\in\N}: X\longrightarrow X$ is a sequence of
contractive maps with the same contraction constant $\al$ and
suppose that $\{f_n\}\longrightarrow f$ (point to point). Then $f$
is $\al$-contractive and for any point $z\in X$ the sequence
$\{f_n\circ\cdots \circ f_1(z)\} \xrightarrow [n\to\infty] \  x_0
$, where $x_0$ is the unique fixed point of $f$.

\noindent Our framework is the following: we consider $\K$  a
field of characteristic zero and the ring of power series
$\K[[x]]$ with coefficients in $\K$. If $g$ is any series  given
by $g=\sum_{n\geq0}g_nx^n$, we recall that  the order of $g$,
$\omega(g)$, is  the smallest nonnegative integer number $n$ such
that $g_n\neq0$ if any exist. Otherwise, that is if  $g=0$, we say
that its order is $\infty$. It is well-known that the space
$(\K[[x]],d)$ is a complete ultrametric space where the distance
between $f$ and $g$ is given by $d(f,g)=\frac{1}{2^{\om(f-g)}}$,
$f,g\in\K[[x]]$. Here we understand that $\frac{1}{2^\infty}=0$.
Moreover the distance between $f$ and $g$ is less than or equal to
$\frac{1}{2^{n+1}}$, i. e.  $d(f,g)\leq\frac{1}{2^{n+1}}$, if and
only if their n-degree Taylor polynomials are equals,
$T_n(f)=T_n(g)$. Finally the sum and product of series are
continuous if we consider the corresponding product topology in
$\K[[x]]\times\K[[x]]$. See for example \cite{Rob00}, \cite{teo}
and \cite{BanPas}  for these topics.

\noindent In this paper $\N$ represents the set of natural numbers
including 0.

\noindent This work is motivated by the following question:
\que\label{Q:question1} Can we sum the arithmetic-geometric series
$\sum_{k=1}^{\infty}kx^{k-1}$ using the Banach Fixed Point Theorem?.

\vspace{0.5cm} \noindent We can sum the geometric series using
{\bf BFPT}. A visual proof of this fact can be found in
\cite{2000Group}. Herein  we recall an analytic proof: (the
peculiar name of the following function will be justified later
on). We consider
\[
\begin{matrix}
h_{m,1}:(\K[[x]],d)&\rightarrow&(\K[[x]],d)\\
 t&\mapsto&xt+1
\end{matrix}
\]
\noindent Since  $h_{m,1}$ is contractive, in fact
$\displaystyle{d(h_{m,1}(t_1),h_{m,1}(t_2))\leq\frac{1}{2}d(t_1,t_2)}$},
we iterate at $t=0$ and we obtain:

\noindent $h_{m,1}(0)=1$\\
\noindent $h_{m,1}^2(0)=x+1$\\
\noindent $h_{m,1}^3(0)=x^2+x+1$\\
\noindent $h_{m,1}^4(0)=x^3+x^2+x+1$\\
\noindent that is,
\[h_{m,1}^{n+1}(0)=\sum_{k=0}^nx^k\]
\noindent As the fixed point of $h_{m,1}$ is the solution of
$xt+1=t$, then $t=\frac{1}{1-x}$ \noindent and from  {\bf BFPT} we
induce
\[
h_{m,1}^{n+1}(0)=\sum_{k=0}^{n}x^k\quad\overrightarrow{n\to\infty}\quad
\frac{1}{1-x}
\]
which is the unique fixed point of the function, in this case,
$h_{m,1}(t)\equiv h_{1}(t)= xt+1$.

\noindent Now it is natural to wonder Question \ref{Q:question1}.

\noindent  We organize the paper in the following way:

\noindent In Section \ref{S:TwoAnswers} we use {\bf BFPT} to a
suitable function related to  Question \ref{Q:question1}. We do
not answer the question by this way but we  find an interesting
arithmetical triangle. Later and using {\bf GBFPT} we answer the
question. Actually we construct the whole  Pascal triangle by this
method, see \cite{teo} and \cite{BanPas}.

\noindent In Section \ref{S:tfggroup} we generalize the method
above, to construct the Pascal triangle, finding so a way to
construct   arithmetical triangles $T(f\mid g)$ for any pair of
series $f$ and $g$ with non null independent terms. Using the
usual product of matrices, we  identify  the well-known Riordan
group, see \cite{teo}.

\noindent In the procedure described above we obtain a new
parametrization of the elements in the Riordan group and so a new
notation  different from the usual ones. In Section 4 we try to
justify the use of our notation alternatively to the usual
notation. Our  method of construction and our notation allow us to
explain easily a way to add and delete columns suitably in a
Riordan array to get another one. For a concrete kind of
triangles, those denoted by $T(1\mid a+bx)$, we can calculate the
inverse only adding adequately new columns to those triangles. In
fact it is an {\it elementary operations} method. We end this
section giving expressions for  the so called $A$ and $Z$
sequences of a Riordan array. These expressions are given in
terms of our notation and related to the inversion in the Riordan
group.

\noindent In Section \ref{S:inversa} we give the main new results
of this paper. We display an algorithm to construct the inverse of
a series and we show the relation with the Lagrange Inversion
Formula. In particular we prove that Banach Fixed Point theorem
gives rise to the Lagrange Inversion Formula.

\section{Two answers: a curious triangle and an iterative
method}\label{S:TwoAnswers}

\noindent  To answer the previous question, in this  section we
are going to recall lightly some examples and tools widely studied
in \cite{teo} and \cite{BanPas}.

\noindent It can be easily shown, see page 2257 in \cite{BanPas},
that there are not a first-degree polynomial with coefficients in
$\K[[x]]$, $f(t)=g(x)t+h(x)$, and any point, $x_0$, such that its
$n+1$-iteration coincides with the partial sum of the
arithmetic-geometric series, that is: $\nexists\, f, x_0 \, /
\,\forall n, \, f^{n+1}(x_0)=\sum_{k=0}^{n}(k+1)x^{k}$.

\noindent In view of this, we are going to iterate a polynomial
whose fixed point is the sum of the arithmetic-geometric series,
that is $\sum_{k=1}^{\infty}kx^{k-1}=\frac{1}{(1-x)^2}$. Since the
equality $t=\frac{1}{(1-x)^2}$ can be converted to
$t=1+(2x-x^2)t$, we consider the polynomial $f(t)=1+(2x-x^2)t$,
with coefficients in $\K[[x]]$,   and we initiate  the iteration
process  at $t=0$:

$f(0)=1$

$f^2(0)=1+2x\mathbf{-x^2}$

$f^3(0)=1+2x+3x^2\mathbf{-4x^3+x^4}$

$f^4(0)=1+2x+3x^2+4x^3\mathbf{-11x^4+6x^5-x^6}$

$f^5(0)=1+2x+3x^2+4x^3+5x^4\mathbf{-26x^5+23x^6-8x^7+x^8}$

$f^6(0)=1+2x+3x^2+4x^3+5x^4+6x^5\mathbf{-57x^6+72x^7-39x^8+10x^{9}-x^{10}}$

\noindent We know that the sequence of iterations converges to the
sum of the arithmetic-geometric series, because $f$ is contractive
in $(\K[[x]], d)$, but we can observe that in each iteration the
partial sum of such series appears plus a remainder. We want to
control the difference with the partial sum. To do this, we display
the coefficients of the remainder as a matrix, that is:
\[
\begin{pmatrix}-1& && &&& & \\ -4&1&& &&&&\\ -11&6&-1&&&&&\\
-26&23&-8&1&&&&\\ -57&72&-39&10&-1&&&\\
-120&201&-150&59&-12&1&&\\
\vdots&\vdots&\vdots&\vdots&\vdots&\vdots&\ddots&\\
\end{pmatrix}
\]

\noindent Observing this matrix we recognize:

\noindent (1) The  rule of construction is similar to that of Pascal
triangle: each element is twice the above element minus the element
above to the left side, that is: $a_{n,k}=2a_{n-1,k}-a_{n-1,k-1}$.

\noindent(2) The elements in the first column are Eulerian numbers
except for the sign.

\noindent(3) The sum of the elements in any row are triangular
numbers with negative sign.

\noindent(4) For every element, the sum of all elements in its row
to the right  and all elements above in its column is zero. That
is, $\sum_{k=1}^{i-1}a_{kj}+\sum_{k=j+1}^{n}a_{ik}=0$.

\noindent(5) The general term is
$a_{n,j}=n+j-1+\sum_{k=1}^{j-1}(-1)^k\binom{n+j-1-k}{n+j-2k}2^{n+j-2k}$.
Etc.

\noindent For an exhaustive development of this triangle see
\cite{BanPas}.

\noindent The above approach, using  {\bf BFPT}, does not give us an
exact answer to our question. To find an adequate answer we consider
the {\bf GBFPT}. This is our way to do this:

\noindent For computability facts we consider the sequence of
functions, with polynomial coefficients, given by
$h_{0,2}(t)=xt$, $h_{1,2}(t)=xt+x$, $h_{2,2}(t)=xt+x+x^2$,
$h_{3,2}(t)=xt+x+x^2+x^3$,
\[
h_{m,2}(t)=xt+x\sum_{k=0}^{m-1}x^k
\]

\noindent Each function $h_{m,2}$ is $\frac{1}{2}$-contractive, so
$ \{h_{m,2}\}$ is an equi-contractive sequence of one-degree
polynomials that converges to $h_2(t)=xt+\frac{x}{1-x}$, i. e. :
\[ \{h_{m,2}\}\longrightarrow h_2(t)=xt+\frac{x}{1-x}
\]
\noindent It is easy to see that  the crossed iterations induced by
{\bf GBFPT} are just the  corresponding partial sums of the
arithmetic-geometric series:

\noindent $h_{0,2}(0)=0$\\
$h_{1,2}(h_{0,2}(0))=x$\\
$h_{2,2}(h_{1,2}(h_{0,2}(0)))=x+2x^2$\\
$h_{3,2}(h_{2,2}(h_{1,2}(h_{0,2}(0))))=x+2x^2+3x^3$\\ \noindent Now
using again {\bf GBFPT} we obtain that,  since  $xt+\frac{x}{1-x}=t$
$\Rightarrow$ $t=\frac{x}{(1-x)^2}$, then
\[
(h_{m,2}\circ \cdots \circ h_{0,2})(0)\longrightarrow
\frac{x}{(1-x)^2}
\]
\noindent these crossed iterations at zero converge to the unique
fixed point of $h_2$, that is, the sum of the arithmetic-geometric
series. So the answer to our question is yes if we are allowed to
use the generalized version of the {\bf BFPT}.

\noindent Recall that  the Pascal triangle is given by:
\[
 \begin{smallmatrix}
 1&& & &&&&& &
 \\
  1&1& & &&&& &&
   \\
  1&2&1& && && && 
  \\
  1&3&3&1&&
&&&&
\\
 1&4&6&4&1&&&&
 \\
1&5&10&10&5&1&&&&
\\ 1&6&15&20&15&6&1&&&
\\
\vdots&\vdots&\vdots&\vdots&\vdots&\vdots&\vdots&
\ddots&\\
\binom{n}{0}&\binom{n}{1}&\binom{n}{2}&\binom{n}{3}&\binom{n}{4}&\binom{n}{5}&
\binom{n}{6}
&\cdots&\binom{n}{n}&
\\
\vdots&\vdots&\vdots&\vdots&\vdots&\vdots&\vdots&\vdots&\vdots&&\\
\frac{1}{1-x}&\frac{x}{(1-x)^2}&\frac{x^2}{(1-x)^3}&\frac{x^3}{(1-x)^4}&\frac{x^4}{(1-x)^5}&
\frac{x^5}{(1-x)^6}&\frac{x^6}{(1-x)^7}
&\cdots&\frac{x^{n-1}}{(1-x)^n}&
\end{smallmatrix}
\]
\noindent We have  just constructed the first two columns of Pascal
tringle using  {\bf BFPT}. In fact we needed only {\bf BFPT} to
construct the first one and  {\bf GBFPT} to get  the second one. The
main observation is that we can follow this iterative procedure to
construct all columns.  For   example we can repeat the process to
construct the third column. To get this goal, we interpret the above
equicontractive sequence $h_{m,2}$ in the following way \[
h_{m,2}(t)=xt+x\sum_{k=0}^{m-1}x^k=xt+xT_{m-1,1}\] \noindent where
$T_{m-1,1}$ is the $m-1$ degree Taylor polynomial of the first
column (which is the geometric series). So in a similar way we
consider the following equicontractive sequence:
\[h_{m,3}(t)=xt+x\sum_{k=0}^{m-1}kx^k=xt+xT_{m-1,2}\]
\noindent where $T_{m-1,2}$ is the Taylor polynomial of the second
column (which  is the arithmetic-geometric series). So, as one can
easily prove,  the crossed iterations for  this sequence coincide
with the partial sums of the third column:

\noindent $h_{0,3}(0)=0$\\
$h_{1,3}(h_{0,3}(0))=0$,\\
 $h_{2,3}(h_{1,3}(h_{0,3}(0)))=x^2$, \\
 $h_{3,3}(h_{2,3}(h_{1,3}(h_{0,3}(0))))=x^2+3x^3$, \\
  $h_{4,3}(h_{3,3}(h_{2,3}(h_{1,3}(h_{0,3}(0)))))=x^2+3x^3+6x^4$, \\
\noindent  Using once more {\bf GBFPT}, we obtain that these crossed
iterations converge to the unique fixed point of the limit function
$h_3(t)=xt+x\frac{x}{(1-x)^2}$. Since
\[h_3(t)=xt+x\frac{x}{(1-x)^2}\ \Rightarrow \
t=\frac{x^2}{(1-x)^3}\] then
\[(h_{m,3}\cdots
h_{0,3})(0)=\sum_{k=0}^{m}\binom{k}{2}x^{k}\longrightarrow\frac{x^2}{(1-x)^3}\]

\noindent Actually, as we said before,  we can construct every
column of Pascal triangle using this process:

\begin{prop} For $n\ge2$, the n-column in Pascal's triangle is obtained
from the (n-1)-column applying the crossed iterations in {\bf GBFPT}
to the sequence $\{h_{m,n}\}_{m\in\N}$ where
\[h_{m,n}(t)=xt+xT_{m-1,n-1}\]  being $T_{m-1,n-1}$ the
(m-1)-Taylor polynomial of the (n-1)-column.
\end{prop}

\section{THE GROUP OF ALL ARITHMETICAL TRIANGLES $T(f\mid g)$.}\label{S:tfggroup}

\noindent Now we generalize the previous iterative method for any
pair of series $f=\sum_{n\geq0}f_nx^n$ and $g=\sum_{n\geq0}g_nx^n$
such that $f_0\neq0$ and $g_0\neq0$. We construct the following
arithmetical triangle $T(f\mid g)$, in this notation the Pascal
triangle is $T(1\mid 1-x)$, where the role of series $1$ is played
by the series $f$ and the role of $1-x$ by $g$.

\vspace{0.5cm}

\centerline{\begin{tabular}{l|llll}
$f_0$ & & & & \\ $ f_1$&$d_{0,0}$& & & \\
$f_2$&$d_{1,0}$&$d_{0,1}$& & \\
$f_3$&$d_{2,0}$&$d_{2,1}$ &$d_{2,2}$&\\
$\vdots$&$\vdots$&$\vdots$&$\vdots$&$\ddots$\\
$f$&$\frac{f}{g}$&$\frac{xf}{g^2}$&$\frac{x^2f}{g^3}$&$\cdots$\\
\end{tabular}}

\vspace{0.5cm}

\noindent In \cite{teo} we interpreted the calculation of
 $\displaystyle{\frac{f}{g}}$  as a fixed point problem. Consider
the sequence
\[
h_{m,1}(t)=T_m\left(\frac{g_0-g}{g_0}\right)t+T_m\left(\frac{f}{g_0}\right)
\]
\noindent where $T_m(f)$ is the m degree Taylor polynomial of $f$.
Observe that the sequence of  crossed iterations has as its limit
the unique fixed point of
$\displaystyle{h_{1}(t)=\frac{g_0-g}{g_0}t+\frac{f}{g_0}}$, that
is: $\displaystyle{\frac{f}{g}}$.
It is the first column of $T(f\mid g)$.

\noindent To construct the second column and the next ones  we
consider the  equicontractive sequences
\[
h_{m,n}(t)=T_m\left(\frac{g_0-g}{g_0}\right)t+xT_{m-1}\left(\frac{x^{n-2}f}{g_0g^{n-1}}\right)
\]
\noindent their corresponding limits are
$h_n(t)=\left(\frac{g_0-g}{g_0}\right)t+x\left(\frac{x^{n-2}f}{g_0g^{n-1}}\right)$
whose corresponding unique fixed points are
\[ t_{n}=\frac{x^{n-1}f}{g^{n}}\]
\noindent The series $t_{n}$ is just  the n-column of our $T(f\mid
g)$.

\begin{thm} Let $\displaystyle{f=\sum_{n\geq0}f_nx^n}$ and
$\displaystyle{g=\sum_{n\geq0}g_nx^n}$ with $g_0\neq0$, then the
Riordan matrix $T(f\mid g)=(d_{n,k})$ is given by

if $k=0$, then
\[d_{0,0}=\frac{f_0}{g_0},
\qquad
d_{n,0}=-\frac{g_1}{g_0}d_{n-1,0}-\frac{g_2}{g_0}d_{n-2,0}\cdots-\frac{g_{n}}{g_0}d_{0,0}+\frac{f_{n}}{g_0}
\]
if $k>0$, then
\[
d_{n,k}=-\frac{g_1}{g_0}d_{n-1,k}-\frac{g_2}{g_0}d_{n-2,k}\cdots-\frac{g_{n-k}}{g_0}d_{k,k}+\frac{d_{n-1,k-1}}{g_0}
\]
\end{thm}

\noindent This theorem gives us the following algorithm to
construct for columns any Riordan matrix:

\begin{alg} Given $\displaystyle{f=\sum_{n\geq0}f_nx^n}$ and
$\displaystyle{g=\sum_{n\geq0}g_nx^n}$ with $g_0\neq0$:

{\bf Step 1:} Calculate the first column $d_{n,0}$.
\[d_{0,0}=\frac{f_0}{g_0},
\qquad
d_{n,0}=-\frac{g_1}{g_0}d_{n-1,0}-\frac{g_2}{g_0}d_{n-2,0}\cdots-\frac{g_{n}}{g_0}d_{0,0}+\frac{f_{n}}{g_0}
\]

{\bf Step k:} Calculate the $k$-column, $d_{n,k}$, using
$k-1$-column.
\[
d_{n,k}=-\frac{g_1}{g_0}d_{n-1,k}-\frac{g_2}{g_0}d_{n-2,k}\cdots-\frac{g_{n-k}}{g_0}d_{k,k}+\frac{d_{n-1,k-1}}{g_0}
\]
\end{alg}
\noindent We can write this algorithm in an informal pseudo-code:

\noindent \textit{READ (f,g,n)\\
SET (d,aux)\\
CALCULATE d[0,0]=f[0]/g[0]\\
$\%$ We calculate the first column\\
FOR i=1 to n\\
FOR k=1 to n\\
CALCULATE  aux(k,i)=g[i-k]*d[k,0]\\
END\\
CALCULATE  d(i,0)=1/g[0]*(f[i]-SUM(aux(:,i)))\\
END\\
$\%$ We calculate the remaining columns\\
FOR j=1 to n\\
FOR i=1 to n\\
FOR k=1 to i\\
CALCULATE  aux(k,i)=g[i-k]*d[k,j]\\
END\\
CALCULATE  d(i,j)=1/g[0]*(d(i-1,j-1)-SUM(aux(:,i))\\
END\\
END\\
PRINT(f,g,d)}
%
%

\noindent So to  construct the  arithmetical triangle $T(f\mid g)$
it is enough to know the ordered pair of series $f$ and $g$, i. e.
the data,  and the algorithm of dividing  two series. Every column
is constructed by the same rule as that in
$\displaystyle{\frac{f}{g}}$ but the coefficients of
$\displaystyle{\frac{f}{g}}$ are replaced with the coefficients of
the previous column. Except for the first column, here we need and
auxiliary column, the coefficients of $f$.

\noindent We can consider the  matrix $T(f\mid g)$, like in Linear
Algebra, as the associated matrix to a $\K$-linear continuous
function, see \cite{teo}:
\[
\begin{matrix}
T(f \mid g):&(\K[[x]],d)&\rightarrow&(\K[[x]],d)\\
&h&\mapsto&T(f\mid g)(h)=\frac{f}{g}h\left(\frac{x}{g}\right)
\end{matrix}
\]
\noindent Using the classical definition of composition of maps and
the behavior of the associated matrix, we can easily find the
formulas for the product and the inverse for these triangles.
\[T(f_1\mid g_1)T(f_2\mid g_2)=T\left(f_1f_2\left(\frac{x}{g_1}\right)\Big|
g_1g_2\left(\frac{x}{g_1}\right)\right)\]

\[(T(f\mid g))^{-1}=T\left(\frac{1}{f\left(\om^{-1}\right)}\Big|
\frac{1}{g\left(\om^{-1}\right)}\right),
\quad\om=\left(\frac{x}{g}\right), \quad
 \om\circ
\om^{-1}= \om^{-1}\circ\om=x
\]

\noindent So if we consider the set of the all arithmetical
triangles with $f_0\neq0$ and $g_0\neq0$ and the usual product of
matrices we obtain a group. Actually this group is the well-known
Riordan group.

\section{On the $T(f\mid g)$ notation.}\label{S:tfgnotation}

\noindent We have received some critics about our notation.
Someone could think that our notation is, in some sense,
cumbersome. Of course it depends strongly on the way you approach
or you run into this group. In this section we are going to give
some reasons why our notation could be very adequate. The basic
formula relating our to the classical notation is
\[
(d(x),h(x))=T\left(\frac{xd}{h}\Big|\frac{x}{h}\right)=\left(\frac{f(x)}{g(x)},\frac{x}{g}\right)=T(f\mid
g)=(d_{i,j})_{i,j\geq0}
\]

\noindent The fundamental equality with our notation is:
\begin{equation}\label{E:FundamentalEquality}
    T(f|g)=T(f|1)T(1|g)
\end{equation}
 \noindent This equality in the other notation is:
\[\left(\frac{f(x)}{g(x)},\frac{x}{g(x)}\right)=(f(x),x)\left(\frac{1}{g(x)},\frac{x}{g(x)}\right)\]
\noindent In (\ref{E:FundamentalEquality}) we can see that every
element of the Riordan group \cite{Sha91} can be expressed by
means of the product of a lower triangular Toepliz matrix whose
columns are the coefficients of series $f$, shifted conveniently,
the matrix $T(f\mid 1)$, and a renewal array, the matrix $T(1\mid
g)$ described by Rogers in \cite{Rog78}. These last kind of
matrices are really similar to the Jabotinsky matrices, see
\cite{Lavoie}. We want to point out that the structure of every
element of the Riordan group is \textit{essentially} in the
structure of the matrix $T(1\mid g)$. For example, to know a
closed formula for the general term of $T(1\mid g)$ gives us at
once a closed formula for the general term in $T(f\mid g)$. This
matrix $T(1\mid g)$, for us, is intrinsically related to the
calculation of $\displaystyle{\frac{1}{g}}$, which is its first
column.

\noindent A comparative table of both  notations is given:


{\tiny

\begin{center}

\begin{tabular}{| c | c | c |}   \hline

{\bf Name} & {\bf(d(t),th(t))}& \textbf{T(f}$\mid$\textbf{g)}\\

\hline

Identity & $\displaystyle{\left(1,t\right)}$ & $T(1\mid 1)$ \\

\hline

Pascal & $\displaystyle{\left(\frac{1}{1-t},\frac{t}{1-t}\right)}$ & $T(1\mid 1-t)$ \\

%
%

%

\hline

Appel subgroup element & (d(t),t) & $T(d\mid 1)$\\
\hline

Associated subgroup element & (1,th(t)) & $\displaystyle{T\left(\frac{1}{h}\mid \frac{1}{h}\right)}$\\
\hline

Bell subgroup element & (d(t),td(t)) & $\displaystyle{T\left(1\mid \frac{1}{d}\right)}$\\
\hline
%
%
%
%

%


\end{tabular}

\end{center}

}

\vspace{0.5cm}

\noindent A curious and symbolically important, for us,  property
of our notation is the way to give, by means of the parameters,
the natural powers of the Pascal triangle:
\begin{prop} For every $n\in\N$, we have
$ T^n(1\mid 1-x)=T(1\mid 1-nx)$
\end{prop}
\begin{proof}
Let us proceed by induction.

For $n=2$: $T^2(1\mid 1-x)=T(1\mid 1-x)T(1\mid 1-x)$,
$\om=\frac{x}{1-x}$. So \[T^2(1\mid 1-x)=T(1\mid
(1-x)(1-\frac{x}{1-x})=T(1\mid 1-2x)\] Suppose that $
T^{n-1}(1\mid 1-x)=T(1\mid 1-(n-1)x)$, then \[ T^n(1\mid
1-x)=T(1\mid 1-x)T^{n-1}(1\mid 1-x)=T(1\mid 1-x)T(1\mid
1-(n-1)x)=\]\[=T(1\mid (1-x)(1-(n-1)\frac{x}{1-x}))=T(1\mid
1-nx)\]
\end{proof}

\vspace{0.5cm}

\noindent Another reason why for us our notation is natural, is
related to the way we begun to study these topics. One of the
first things we did was to find our curious triangle described in
Section \ref{S:TwoAnswers}. From our notation, the description as
a Riordan array is:
\[
\begin{pmatrix}-1& && &&& & \\ -4&1&& &&&&\\ -11&6&-1&&&&&\\
-26&23&-8&1&&&&\\ -57&72&-39&10&-1&&&\\
-120&201&-150&59&-12&1&&\\
\vdots&\vdots&\vdots&\vdots&\vdots&\vdots&\ddots&\\
\end{pmatrix}=T\left(\frac{1}{(1-x)^2}\mid2x-1\right)
\]
\noindent This notation resembles both the problem we were treating
and the algorithm of construction.

\noindent Another thing we can describe easily with our notation
is the fact that,  with our construction method by columns, we can
add  new columns to the left for every element of the Riordan
group to obtain again a Riordan array intrinsically related to the
initial one, for example:
\[
\begin{pmatrix}1&& && &&& & \\2&-1& && &&& & \\3& -4&1&& &&&&\\4& -11&6&-1&&&&&\\
5&-26&23&-8&1&&&&\\6& -57&72&-39&10&-1&&&\\
7&-120&201&-150&59&-12&1&&\\
\vdots&\vdots&\vdots&\vdots&\vdots&\vdots&\vdots&\ddots&\\
\end{pmatrix}=T\left(\frac{2x-1}{(1-x)^2}\mid2x-1\right)
\]

\noindent Note that we added a new column to the left and look at
the way the parameters changed in our notation.

\noindent In general we can construct a family of new Riordan
matrices closely related to it. For example by definition of
Riordan array we get
\[
T(fg\mid g)=\left(%
\begin{array}{cccccc}
 f_0&  &  &  &  &    \\
 f_1& d_{0,0} &  &   &  &  \\
 f_2& d_{1,0} & d_{1,1}  &  &  &  \\
 f_3& d_{2,0} & d_{2,1} & d_{2,2}  &  &  \\
 f_4& d_{3,0} & d_{3,1} & d_{3,2} & d_{3,3}  &  \\
  \vdots & \vdots & \vdots & \vdots & \vdots & \ddots   \\
\end{array}%
\right)
\]

\[
T\left(\frac{f}{g}\mid g\right)=\left(%
\begin{array}{cccccc}
  d_{1,1} &  &  &  & & \\
  d_{2,1} & d_{2,2} &  &  &  &\\
  d_{3,1} & d_{3,2} & d_{3,3} &  & & \\
  d_{4,1} & d_{4,2} & d_{4,3} & d_{4,4} &  \\
  d_{5,1} & d_{5,2} & d_{5,3} & d_{5,4} & d_{5,5}  \\
  \vdots & \vdots & \vdots & \vdots & \vdots & \ddots \\
\end{array}%
\right)
\]
where $\displaystyle{f=\sum_{n\geq0}f_nx^n}$ and $T(f\mid
g)=(d_{n,k})_{n,k\in\N}$.

\noindent By the same way we observe that
$T\left(\frac{f}{g^m}\mid g\right)$ for $m\in\N$ is the matrix
obtained from $T(f\mid g)$ by deleting the first $m$-rows and
$m$-columns. Moreover $T(fg^m\mid g)$ is the unique Riordan matrix
with the property that by deleting the first $m$-rows and
$m$-columns from it, we obtain $T(f\mid g)$. In fact it can be
easily proved the following:

\begin{prop}
Let $T(f\mid g)=(d_{n,k})_{n,k\in\N}$ be a Riordan matrix, and
$m\in\Z$ then $T(fg^m\mid g)=(\tilde{d}_{n,k})_{n,k\in\N}$, with
$\tilde{d}_{n,k}=[x^{n-k}]fg^{m-k-1}$. Where $[x^j]S$ stands for
the $j$-coefficient of the formal power series $S$.
\end{prop}

We can always embed any $T(f|g)$ in a bi-infinite lower triangular
matrix.

For our example we have
\[
\left(
  \begin{array}{cccccccccccc}
    &   & & & & & \vdots  &  & &  &  \\
   &  -1 & 0 & 0 & 0 & 0 & 0 & 0 & 0 & 0 & 0 \\
  &   8 & 1 & 0 & 0 & 0 & 0 & 0 & 0 & 0 & 0 \\
  &   -23 & -6 & -1 & 0 & 0 & 0 & 0 & 0 & 0 & 0 \\
  &   26 & 11 & 4 & 1 & 0 & 0 & 0 & 0 & 0 & 0 \\
  \cdots &  -5 & -4 & -3 & -2 & -1 & 0 & 0 & 0 & 0 & 0 &\cdots\\
   & -4 & -3 & -2 & -1 & 0 & 1 & 0 & 0 & 0 & 0\\
  &   -3 & -2 & -1 & 0 & 1 & 2 & -1 & 0 & 0 & 0 \\
  &   -2 & -1 & 0 & 1 & 2 & 3 & -4 & 1 & 0 & 0 \\
  &   -1 & 0 & 1 & 2 & 3 & 4 & -11 & 6 & -1 & 0 \\
  &   0 & 1 & 2 & 3 & 4 & 5 & -26 & 23 & -8 & 1 \\
    &   & & & & &\vdots   &  & &  &  \\
  \end{array}
\right)
\]
%

\noindent This construction is very nice in  the case of Riordan
matrices of the kind $T(1\mid a+bx)$, (treated as change of
variables in \cite{BanPas}). For these matrices we have
\[
\left(
  \begin{array}{ccccccccccc}
  &&  \vdots & \vdots & \vdots & \vdots & \vdots & \vdots & \vdots &  \\
  \frac{1}{a}\left(\frac{a}{1-bx}\right)^4 &\leftarrow&  a^3 & 0 & 0 & 0 & 0 & 0 & 0 & \cdots \\
   \frac{1}{a}\left(\frac{a}{1-bx}\right)^3 &\leftarrow &  3a^2b& a^2 & 0 & 0 & 0 & 0 & 0 & \cdots\\
    \frac{1}{a}\left(\frac{a}{1-bx}\right)^2&\leftarrow&  3ab^2 & 2ab & a & 0 & 0 & 0 & 0 & \cdots \\
    \frac{1}{a}\frac{a}{1-bx}&\leftarrow &  b^3 & b^2 & b & 1 & 0 & 0 & 0 & \cdots\\
   \frac{1}{a}&\leftarrow &  0 & 0 & 0 & 0  & \frac{1}{a} & 0 & 0& \cdots \\
   \frac{1}{a}\left(\frac{1-bx}{a}\right) &\leftarrow &  0 & 0 & 0 & 0  & -\frac{b}{a^2} & \frac{1}{a^2} & 0& \cdots \\
   \frac{1}{a}\left(\frac{1-bx}{a}\right)^2&\leftarrow &  0 & 0 & 0 & 0 & \frac{b^2}{a^3} & -\frac{2b}{a^3} & \frac{1}{a^3} & \cdots \\
   & &  \downarrow & \downarrow  & \downarrow  & \downarrow  & \downarrow  & \downarrow  & \downarrow  & \vdots    \\
   &&  (a+bx)^3 & (a+bx)^2 & (a+bx) & 1 & \frac{1}{a+bx} & \frac{1}{(a+bx)^2} & \frac{1}{(a+bx)^3} & \cdots \\
  \end{array}
\right)
\]

\noindent Note that, as in the Pascal triangle, we can see the
above matrix in the following way
\[\left(
              \begin{array}{cc}
               \frac{1}{a} T^{-1}(1\mid a+bx)\blacklozenge & 0 \\
                0 & T(1\mid a+bx) \\
              \end{array}
            \right)
\]
\noindent Where $\frac{1}{a} T^{-1}(1\mid a+bx)\blacklozenge $ is
$\frac{1}{a} T^{-1}(1\mid a+bx)$ placed in the same way as the
inverse of the Pascal triangle is placed in the Hexagon of Pascal in
page 194 in \cite{H-P2}. Note that $T^{-1}(1\mid
a+bx)=T(1\mid\frac{1-bx}{a})$. This gives us a method to calculate
$T^{-1}(1\mid a+bx)$ by means of elementary operations.

\vspace{0.5cm}

\noindent Our algorithm of construction do not need the so called
$A$ and $Z$ sequences of a Riordan array, see \cite{Merlini} and
\cite{Spr94}. But, in our notation, they appear in the expression
of the inverse of the Riordan array $T(f\mid g)$ giving us
specially aesthetic formula:

\begin{prop} Let $\displaystyle{f=\sum_{n\geq0}f_nx^n}$ and
$\displaystyle{g=\sum_{n\geq0}g_nx^n}$ be two formal power series
with $f_0\neq0$ and $g_0\neq0$. Suppose that $A$ and $Z$ represent
the $A$-sequence and the $Z$-sequence, respectively, of $T(f\mid
g)$. Then
\begin{itemize}
  \item [(i)] $T^{-1}(1\mid g)=T(1\mid A)$
  \item [(ii)]$T^{-1}(f\mid g)=T\left(\frac{g_0}{f_0}(A-xZ)\big|
  A\right)$
\end{itemize}
\end{prop}
\begin{proof} (i) From Theorem 1.3 of \cite{Spr94}, the $A$-sequence is the unique
series, with $A(0)\neq0$,  such that
$\displaystyle{\frac{1}{g}=A\left(\frac{x}{g}\right)}$. So
$\displaystyle{\frac{x}{g}=xA\left(\frac{x}{g}\right)}$. If
$\displaystyle{\om=\frac{x}{g}}$ then
$\displaystyle{\om=xA\left(\om\right)}$. On the other hand as
$\om^{-1}\circ \om=x$ and $\displaystyle{\om=\frac{x}{g}}$ then
$\displaystyle{x=\om g}$, composing with $\om^{-1}$ we get
$\displaystyle{\om^{-1}=xg(\om^{-1})}$ and
$\displaystyle{\frac{\om^{-1}}{x}=g(\om^{-1})}$. So composing with
$\om^{-1}$ but now in $\displaystyle{\om=xA\left(\om\right)}$ we
get $\displaystyle{x=\om^{-1}A(x)}$ then
$\displaystyle{1=\frac{\om^{-1}}{x}A(x)}$ so $1=g(\om^{-1})A(x)$
then $\displaystyle{A(x)=\frac{1}{g(\om^{-1})}}$. Since
$T^{-1}(1\mid g)=T(1\mid\frac{1}{g(\om^{-1})})=T(1\mid A)$.

\noindent (ii) From Theorem 2.3 in \cite{Merlini} we obtain that
the $Z$ is determined by the equality

$\displaystyle{\om^{-1}Z=
1-\frac{f_0g(\om^{-1})}{g_0f(\om^{-1})}}$. From here we get
$\displaystyle{\frac{1}{f(\om^{-1})}=\frac{g_0}{f_0}(A-xZ)}$. So
\[
T^{-1}(f\mid
g)=T\left(\frac{1}{f(\om^{-1})}\big|\frac{1}{g({\om^{-1})}}\right)=T\left(\frac{g_0}{f_0}(A-xZ)\big|
  A\right)
\]
\end{proof}

\begin{cor}
\[T^{-1}(f\mid
g)=T(1\mid A)T\left(\frac{1}{f}\mid 1\right)\]
\end{cor}

\section{Lagrange inversion formula via Banach fixed point Theorem}\label{S:inversa}

\noindent In the previous section we showed that to calculate the
inverse of $T(1\mid g)$ we need, in particular, to calculate
$\om^{-1}$, where $\om=\frac{x}{g}$ and then
$\om^{-1}=xg(\om^{-1})$. So we consider the function
$F:x\K[[x]]\rightarrow x\K[[x]]$ defined by $F(y)=xg(y)$. Here
$x\K[[x]]$ represents the series with null independent term. This
function is $\frac{1}{2}$-contractive since
\[
d(F(y_1),F(y_2))=\frac{1}{2^{\omega(xg(y_1)-xg(y_2))}}=\frac{1}{2^{\omega(g(y_1)-g(y_2))+1}}\leq
\frac{1}{2}d(y_1,y_2)
\]
\noindent The domain, $x\K[[x]]$, of $F$ is the closed ball, in
$(\K[[x]],d)$, whose center is the series $S=0$ and the ratio is
$\frac{1}{2}$. Consequently our domain is also complete with the
relative metric. So the unique fixed point of $F$ is
$\om^{-1}=\left(\frac{x}{g}\right)^{-1}$ and {\bf BFPT} can be
applied.

\noindent The {\bf BFPT} gives us a theoretical iterative process
to calculate $\om^{-1}$. To convert this method into an effective
approximation process we first note that the relation
$d(S_1,S_2)\leq\frac{1}{2^{m+1}}$ means that the m degree Taylor
polynomials of both series are equals, that is
$T_m(S_1)=T_m(S_2)$. Then we obtain the following algorithm:

\noindent Suppose $g=\sum_{n\geq0}g_nx^n$. We begin to iterate at
$S=0$. $F(0)=g_0x$. This means that $T_1(\om^{-1})=g_0x$. Using
the $\frac{1}{2}$-contractivity of $F$ we get
$T_2(F(g_0x))=T_2(\om^{-1})$. Since
$F(g_0x)=g_0x+g_0g_1x^2+\cdots$ we obtain
$T_2(\om^{-1})=g_0x+g_0g_1x^2$. Similar arguments allow us to
prove that $T_3(F(g_0x+g_0g_1x^2))=T_3(\om^{-1})$. The above
construction can be summarized in the following (the notation is
as above):

\begin{prop}\label{A:inversa}
\[T_m(\om^{-1})=T_m(F(T_{m-1}(F(\cdots(F(T_1(F(0))))\cdots))))\]
\end{prop}

\noindent Following this process we get

\noindent $T_1(\om^{-1})=g_0x$\\
\noindent $T_2(\om^{-1})=g_0x+g_0g_1x^2$\\
\noindent $T_3(\om^{-1})=g_0x+g_0g_1x^2+(g_0g_1^2+g_0^2g_2)x^3$\\
\noindent $T_4(\om^{-1})=g_0x+g_0g_1x^2+(g_0g_1^2+g_0^2g_2)x^3+(g_0g_1^3+3g_0^2g_1g_2+g_0^3g_3)x^4$\\
\noindent
$T_5(\om^{-1})=g_0x+g_0g_1x^2+(g_0g_1^2+g_0^2g_2)x^3+(g_0g_1^3+3g_0^2g_1g_2+g_0^3g_3)x^4+
(g_0g_1^4+6g_0^2g_1^2g_2+2g_0^3g_2^2+4g_0^3g_1g_3+g_0^4g_4)x^5$

\noindent  If we recall the Cauchy powers of the series $g$:

\noindent $g(x)=\mathbf{g_0}+g_1x+g_2x^2+g_3x^3+g_4x^4+\cdots$\\
\noindent
$g^2(x)=g_0^2+\mathbf{2g_0g_1}x+(2g_0g_2+g_1^2)x^2+(2g_0g_3+2g_1g_2)x^3\cdots$\\
\noindent $g^3(x)=g_0^3+3g_0^2g_1x+\mathbf{3(g_0g_1^2+g_0^2g_2)}x^2+(6g_0g_1g_2+3g_0^2g_3+g_1^3)x^3+\cdots$\\
\noindent $g^4(x)=g_0^4+4g_0^3g_1x+(4g_0^3g_2+6g_0^2g_1^2)x^2+\mathbf{4(g_0g_1^3+3g_0^2g_1g_2+g_0^3g_3)}x^3+\cdots$\\
\noindent
$g^5(x)=g_0^5+5g_0^4g_1x+(5g_0^4g_2+10g_0^3g_1^2)x^2+(20g_0^3g_1g_2+10g_0^2g_1^3+5g_0^4g_3)x^3+
\mathbf{5(g_0g_1^4+6g_0^2g_1^2g_2+2g_0^3g_2^2}$

\noindent $\mathbf{+4g_0^3g_1g_3+g_0^4g_4})x^4+\cdots$

\noindent comparing adequately the coefficients of $\om^{-1}$ and
the powers of $g$ we obtain the next relationships:

\noindent $[x]\om^{-1}=[x^0]g$\\
\noindent $[x^2]\om^{-1}=\frac{1}{2}[x^1]g^2$\\
\noindent $[x^3]\om^{-1}=\frac{1}{3}[x^2]g^3$\\
\noindent $[x^4]\om^{-1}=\frac{1}{4}[x^3]g^4$\\
\noindent $[x^5]\om^{-1}=\frac{1}{5}[x^4]g^5$\\

\noindent These equalities allow us to predict and motivate the
classical Lagrange Inversion Formula, see \cite{Stanley} page 36:
\[
[x^{n+1}]\om^{-1}=\frac{1}{n+1}[x^n]g^{n+1}, \ \text{with} \
\om=\frac{x}{g}
\]
\noindent  From now on we denote $T_j\equiv T_j(\om^{-1})$. To
show how this process works note that
\[
F(T_n)=x(g_0+g_1T_n+g_2T_n^2+\cdots+g_nT_n^n+\cdots)=\]
\[
=T_n+(g_1[x^n]T_n+g_2[x^n]T_n^2+\cdots+g_n[x^n]T_n^n)x^{n+1}+S_{n+2}\quad
\text{with}\quad S_{n+2}\in x^{n+2}\K[[x]]
\]
\noindent So
\[
[x^{n+1}]F(T_n)=\sum_{k=1}^n[x^k]g[x^n](\om^{-1})^k
\]
\noindent Because $[x^n](\om^{-1})^k=[x^n](T_n)^k$ for any $k\leq
n$. Suppose now that we know
\[
n[x^n](\om^{-1})^k=k[x^{n-k}]g^n \qquad \text{for}\quad k\leq n
\]
\noindent then
\[
[x^{n+1}]F(T_n)=[x^{n+1}]\om^{-1}=\frac{1}{n}\sum_{k=1}^n
k[x^k]g[x^{n-k}]g^n=\frac{1}{n}[x^{n-1}]g'g^n=\frac{1}{n+1}[x^n]g^{n+1}
\]
\noindent Note that in the above development we need to know
$\displaystyle{n[x^n](\om^{-1})^k=k[x^{n-k}]g^n \ \text{for}\
k\leq n}$. In fact we can give a proof of all above using
essentially the fact that $\om^{-1}$ is a fixed point of certain
contractive function.
\begin{thm}(Lagrange inversion via Banach Fixed Point Theorem)
Let $\K$ be a field of characteristic zero. Suppose that $\om$ is
a formal power series in $\K[[x]]$ with $\om(0)=0$ and
$\om'(0)\neq0$. Then
\[
n[x^n](\om^{-1})^k=k[x^{n-k}]\left(\frac{x}{\om}\right)^n \qquad
\text{for}\qquad n,k\in\N
\]
\end{thm}
\begin{proof}
Let $\displaystyle{g=\frac{x}{\om}}$. So $[x^0]g\neq0$. As proved
before $\om^{-1}$ is the unique fixed point of the
$\frac{1}{2}$-contractive function $F:x\K[[x]]\rightarrow
x\K[[x]]$ defined by $F(y)=xg(y)$. Iterating at $y=0$ we get
\[
[x^1]\om^{-1}=[x^1]F(0)=[x^0]g
\]
If $k>1$, note that $[x^1](\om^{-1})^k=0$  and $[x^{1-k}]g=0$ and
then
\[
[x^1]\om^{-1}=k[x^{1-k}]g
\]
Let proceed by induction on $n$. Suppose that
\[
j[x^j](\om^{-1})^k=k[x^{j-k}]g^j \qquad \text{for} \qquad j\leq n,
\quad k\geq1
\]
\noindent Note that there are actually  only a finite number of
suppositions on $k$. Because if $j<k$, then $[x^{j-k}]g^j=0=
[x^j](\om^{-1})^k$. Then the equality holds trivially.

\noindent Since $\om^{-1}=xg(\om^{-1})$, then for any $k$
$(\om^{-1})^k=x^kg^k(\om^{-1})$. Consequently
\[
[x^{n+1}](\om^{-1})^k=[x^{n+1}]x^kg^k(\om^{-1})=[x^{n+1-k}]g^k(\om^{-1})=\sum_{j=0}^{n+1-k}[x^j]g^k[x^{n+1-k}](\om^{-1})^j
\]
\noindent by the induction hypothesis
\[
[x^{n+1}](\om^{-1})^k=\frac{1}{n+1-k}\sum_{j=0}^{n+1-k}j[x^j]g^k[x^{n+1-k-j}]g^{n+1-k}
\]
\noindent Let us call $h=g^k$
\[
[x^{n+1}](\om^{-1})^k=\frac{1}{n+1-k}\sum_{j=0}^{n+1-k}j[x^j]h[x^{n+1-k-j}]h^{\frac{n+1-k}{k}}=
\frac{1}{n+1-k}\sum_{j=1}^{n+1-k}[x^{j-1}]h'[x^{n+1-k-j}]h^{\frac{n+1-k}{k}}=
\]
\[
=\frac{1}{n+1-k}\sum_{j=0}^{n-k}[x^{j}]h'[x^{n-k-j}]h^{\frac{n+1-k}{k}}=\frac{1}{n+1-k}[x^{n-k}](h'h^{\frac{n+1-k}{k}})
=
\]
\[=\frac{1}{n+1-k}[x^{n-k}]\left(\frac{k}{n+1}h^{\frac{n+1-k}{k}}
\right)'=\frac{k}{n+1}[x^{n+1-k}]h^{\frac{n+1}{k}}=\frac{k}{n+1}[x^{n+1-k}]g^{n+1}
\]
\end{proof}

\noindent The development above gives us the following algorithm
to calculate the coefficients of the compositional inverse of
$\displaystyle{\om=\frac{x}{g}}$:

\begin{alg}
Given $\displaystyle{g=\sum_{j\geq0}g_jx^j}$, with $g_0\neq0$.
Given $F(y)=xg(y)$, $y\in x\K[[x]]$.

{\bf step 1}: (Initial.) $T_1=g_0x$.

{\bf step k (2 to n)}: Calculate the Taylor polynomial of order
$k$ of $F(T_{k-1})$
\end{alg}

\noindent We can write this algorithm in an informal pseudo-code:

\noindent\textit{READ (g,F,n)\\
SET T\\
CALCULATE for i from 2 to n do\\
T[i]:=convert(series(F(T[i-1],x=,i)),polynom);\\
end\\
PRINT T}


 {\bf Acknowledgment:} The author thanks Manuel A. Morón for his
 helpful comments. I also thank the referee for his/her suggestions and comments which
 strongly improved earlier version of this paper. The author was partially supported by the grant
 MICINN-FIS2008-04921-C02-02.

\end{document}